\documentclass[a4paper,10pt]{article}
\usepackage{latexsym}
\usepackage{amsfonts}
\usepackage{amsmath}
\usepackage{amsthm}
\usepackage{amssymb}
\usepackage{bbm}

  \usepackage{paralist}
  \usepackage{graphics} 
  \usepackage{epsfig} 
\usepackage{graphicx}  \usepackage{epstopdf}
\usepackage{paralist}
\renewenvironment{description}[0]{\begin{compactdesc}}{\end{compactdesc}}
\usepackage[affil-it]{authblk}
 \usepackage[colorlinks=true]{hyperref}
\hypersetup{urlcolor=blue, citecolor=red}

\newtheorem{theorem}{Theorem}[section]

\newtheorem{proposition}{Proposition}

\newtheorem{definition}[theorem]{Definition}
\newtheorem{example}[theorem]{Example}
\newtheorem{remark}{Remark}

\newcommand{\reals}{\mathbb{R}}

 \newcommand{\C}[1]{\mathbf{C^{#1}}}
 \renewcommand{\L}[1]{\mathbf{L^{#1}}}
 \newcommand{\Lip}{\mathbf{Lip}}
 \newcommand{\lip}{\mathrm{Lip}}
 \newcommand{\Cc}[1]{\mathbf{C_c^{#1}}}
 
 \renewcommand{\d}{\rm{d}}
 \newcommand{\iw}{\textrm{i}_{w}}
 \newcommand{\caratt}[1]{{\displaystyle\chi_{\strut{\textstyle #1}}}}
 \newcommand{\modulo}[1]{{\left|#1\right|}}
 \newcommand{\norma}[1]{{\left\|#1\right\|}}
 \newcommand{\BV}{\mathbf{BV}}
 
 \newcommand{\sign}{\mathrm{sign}}
 
 \newcommand{\Lloc}[1]{\mathbf{L^{#1}_{loc}}}
 
\begin{document}

\title{Riemann problems with non--local point constraints and capacity drop}

\author[1,2]{Boris Andreianov}

\author[2]{Carlotta Donadello}

\author[2]{Ulrich Razafison}

\author[3]{Massimiliano~D.~Rosini}

\affil[1]{Institut f\"ur Mathematik, Technische Universit\"at Berlin,

Str.~des 17.~Juni 136, 10623 Berlin, Germany}

\affil[2]{Laboratoire de Math\'ematiques CNRS UMR 6623,

 Universit\'e de Franche-Comt\'e,  

16 route de Gray, 25030 Besan\c{c}on Cedex, France}

\affil[3]{
ICM, Uniwersytet Warszawski, 

ul.~Prosta 69, 00838 Warsaw, Poland}

 %\email{boris.andreianov@univ-fcomte.fr}
 %\email{carlotta.donadello@univ-fcomte.fr}
 %\email{ulrich.razafison@univ-fcomte.fr}
 %\email{mrosini@icm.edu.pl}

\maketitle

\begin{abstract}
In the present note we discuss in details the Riemann problem for a one--dimensional hyperbolic conservation law subject to a point constraint. We investigate how the regularity of the constraint operator impacts the well--posedness of the problem, namely in the case, relevant for numerical applications, of a discretized exit capacity. We devote particular attention to the case in which the constraint is given by a non--local operator depending on the solution itself. We provide several explicit examples.

We also give the detailed proof of some results announced in the paper [Andreainov, Donadello, Rosini, \emph{Crowd dynamics and conservation laws with non--local point constraints and capacity drop}], which is devoted to  existence and stability for a more general class of Cauchy problems subject to Lipschitz continuous non--local point constraints.
\end{abstract}



\section{Introduction}\label{sec:intro}

\subsection{Point constraints in traffic modeling}

Traffic modeling is an exciting and fast--developing field of research with plentiful applications to real life. While this subject was initially limited to the description and the management of vehicular traffic, we see  a growing interest nowadays on different applications as crowd dynamics and bio--mathematics. This note is related to an extensive on--going research project concerning  the theoretical and the numerical study of macroscopic models for which the definition of solution involves an artificial limitation of the flux in a finite number of points. From the modeling point of view, this may correspond to a narrow exit in crowd modeling, a toll gate in vehicular traffic, a cell membrane in bio-medical modeling.

In the pioneering paper \cite{ColomboGoatinConstraint},  R.~Colombo and P.~Goatin introduced point constraints in the classical one--dimensional LWR~road traffic model~\cite{LighthillWhitham, Richards}, with the goal to model the presence of obstacles on the road as toll gates and road lights. This model reads as
\begin{subequations}\label{eq:constrianed}
\begin{align}\label{eq:constrianed1}
    \partial_t\rho + \partial_xf(\rho) &= 0 & (t,x) &\in \reals_+\times\reals\\
    \label{eq:constrianed2}
    f\left(\rho(t,0\pm)\right) &\le q(t) & t &\in \reals_+\\
    \label{eq:constrianed3}
    \rho(0,x) &=\rho_0(x) & x &\in \reals,
\end{align}
\end{subequations}
where $\rho = \rho(t,x) \in \left[0,R\right]$ is the unknown (mean) density at time $t \in \reals_+$ of vehicles moving along the road parameterized by the coordinate $x \in \reals$. Then, $R \in \reals_+$ is the maximal road density, $f ~\colon~ [0,R] \to \reals$ is the nonlinearity relating the flux in the direction of increasing $x$ to the density, $q ~\colon~ \reals_+ \to \reals_+$ is a given function prescribing the maximal flow allowed through the point $x=0$, and $\rho_0  ~\colon~ \reals \to [0,R]$ is the initial (mean) density. Finally, $\rho(t,0-)$ denotes the left measure theoretic trace along the constraint implicitly defined by
\begin{align*}
    \lim_{\varepsilon\downarrow0}\frac{1}{\varepsilon} \int_0^{+\infty}  \int_{-\varepsilon}^0 \modulo{\rho(t,x) - \rho(t,0-)} ~\phi(t,x) ~{\d} x ~{\d} t&=0
\end{align*}
for all $\phi \in \Cc\infty(\reals^2;\reals)$. The right measure theoretic trace, $\rho(t,0+)$, is defined analogously.

In the above setting, the authors of \cite{ColomboGoatinConstraint} were able to prove existence and well--posedness of solutions. Further theoretical considerations and numerics associated to this model have been developed in \cite{scontrainte}.

In view of the applications, however, we also need to consider the case in which the evolution of the constraint function $q$ is not given beforehand, but  depends on the solution $\rho$ itself in a neighborhood of $x=0$. In such situation we say that the point constraint is non--local. In this way we obtain crowd and cell membrane dynamics models described by coupled PDE--ODE systems for which  the existence and well--posedness of solutions are not trivial matter.  Nevertheless, this case has a practical interest.  In crowd dynamics, as an example, the experimental observations by E.~Cepolina in \cite{Cepolina2009532} prove that the irrational behavior of pedestrians at bottlenecks ends up by reducing the maximal possible outflow. This phenomenon, called \emph{capacity drop}, is also related to other effects observed in crowd dynamics, such \emph{Faster Is Slower} and the Braess' paradox.

In full generality, we may consider the constraint function $q$ as follows
\begin{equation*}
  q(t) = \mathcal{Q}[\rho](t)~,
\end{equation*}
where $\mathcal{Q} : \C0\left([0,T]; \L1\left(\reals; [0,R]\right)\right)\to \L1\left([0,T]; [0,R]\right)$. The minimal regularity properties to impose on $\mathcal{Q}$ in order to achieve well--posedness of solutions are not known at the moment, and they are the object of one of our current research projects.

\subsection{An example of non--local point constraint}

In the paper \cite{AndreianovDonadelloRosini}, B.~Andreianov, C.~Donadello and M.~D.~Rosini proposed a model which generalizes the one in~\cite{ColomboGoatinConstraint}  and  consists of a Cauchy problem for a one--dimensional hyperbolic conservation law as~\ref{eq:constrianed} subject to a  non--local point constraint of the form
\begin{equation}
    \label{eq:constrianed4}
   q(t) = p\left( \int_{\reals_-} w(x) ~ \rho(t,x) ~{\d} x\right) \quad t \in \reals_+ ~.
\end{equation}
Here $p ~\colon~ \reals_+ \to \reals_+$ prescribes the maximal flow allowed through an exit placed in $x=0$ as a function of the weighted average density of pedestrians, $\rho$, in a left neighborhood of the exit and  $w  ~\colon~ \reals_- \to \reals_+$ is the weight function used to  average the density. The authors of~\cite{AndreianovDonadelloRosini} proved well-posedness of the Cauchy problem in $\L\infty(\reals ; \left[0,R\right])$ for the model under the following assumptions on the regularity of $f$, $w$ and $p$, see Figure~\ref{fig:Korn},

\begin{enumerate}[\textbf{(P1)}]
  \item[\textbf{(F)}]  $f \in \Lip\left( [0,R]; \left[0, +\infty\right[ \right)$, $f(0) = 0 = f(R)$ and there exists $\bar\rho \in \left]0,R\right[$ such that $f'(\rho)~(\bar\rho-\rho)>0$ for a.e.~$\rho \in [0,R]\setminus\{\bar\rho\}$.
  \item[\textbf{(W)}] $w \in \L\infty(\reals_-;\reals_+)$ is an increasing map, $\norma{w}_{\L1(\reals_-;\reals_+)} = 1$ and there exists $\iw >0$ such that $w(x) = 0$ for any $x \le -\iw$.
\item[\textbf{(P1)}] $p$ belongs to $\Lip \left( \left[0,R\right]; \left]0,f(\bar\rho)\right] \right)$ and it is a non-increasing map.
\end{enumerate}

We recall that, in the previous literature, the only macroscopic model able to reproduce the capacity drop at bottlenecks is the CR model introduced by R.~Colombo and M.~D.~Rosini in \cite{ColomboRosini1}. The Riemann solver for the model described in \cite{ColomboRosini1} is fairly intricate; in this note, our main goal is to describe in an exhaustive way the Riemann solver (or, rather, solvers) for the model of \cite{AndreianovDonadelloRosini}. Notice that, in contrast to \cite{ColomboRosini1}, a specific nonclassical Riemann solver has to be used only at the exit point $x=0$ while the simple classical Riemann solver is used elsewhere.

The notion of solution we adopt is a natural extension of the one introduced in \cite{ColomboGoatinConstraint}, for a constrained Cauchy problem of the form~\ref{eq:constrianed}.
\begin{definition}\label{def:entropysol}
    Assume conditions~\textbf{(F)}, \textbf{(W)} and that $p$ is a non--increasing, possibly multivalued, map with values in $\left]0, f(\bar \rho)\right]$. A map $\rho \in \L\infty( \reals_+\times\reals;[0,R]) \cap \C0(\reals_+; \Lloc1(\reals;[0,R]))$ is an entropy weak solution to~\ref{eq:constrianed}, \ref{eq:constrianed4} if the following conditions hold:
    \begin{enumerate}
    \item\label{cond:1b} There exists $q \in \L\infty( \reals_+;[0, f(\bar\rho)])$ such that for every test function $\phi \in \Cc\infty(\reals^2; \reals_+)$ and for every $k \in [0,R]$
    \begin{subequations}\label{eq:entropysol}
    \begin{align}\label{eq:entropysol1}
        &\int_{\reals_+} \int_\reals \left[ \modulo{\rho-k} \partial_t\phi + \sign(\rho-k) \left(f(\rho) - f(k)\right) \partial_x\phi \right] ~{\d} x ~{\d} t\\ \label{eq:entropysol2}
        &+ 2 \int_{\reals_+} \left[1 - \dfrac{q\left( t\right)}{f(\bar\rho)}   \right] f(k) ~\phi(t,0) ~{\d} t\\ \label{eq:entropysol3}
        &+ \int_\reals \modulo{\rho_0(x) - k} ~\phi(0,x) ~{\d} x
        \ge 0~,
    \end{align}
    and
    \begin{align}\label{eq:entropysol4}
        &f\left(\rho(t, 0\pm)\right) \le q\left(t\right) \hbox{ for a.e.~}t \in \reals_+~.
    \end{align}
    \end{subequations}
    \item In addition $q$ is linked to $\rho$ by the relation~\ref{eq:constrianed4}.
    \end{enumerate}
\end{definition}
If $q$ is given \textit{a priori}, then~\ref{eq:entropysol} is the definition of entropy weak solution to problem~\ref{eq:constrianed}. We refer to  Proposition~2.6 in \cite{scontrainte} for a series of equivalent formulations of conditions~\ref{eq:entropysol}.

The next proposition lists the basic properties of a entropy weak solution of~\ref{eq:constrianed}, \ref{eq:constrianed4}, for the case of a single valued $p$, the proof is given in \cite{AndreianovDonadelloRosini}.
\begin{proposition}\label{prop:ws}
    Let $[t \mapsto \rho(t)]$ be an entropy weak solution of~\ref{eq:constrianed}, \ref{eq:constrianed4} in the sense of Definition~\ref{def:entropysol}. Then
    \begin{enumerate}[(1)]
      \item It is also a weak solution of the Cauchy problem~\ref{eq:constrianed1}, \ref{eq:constrianed3}.
      \item Any discontinuity satisfies the Rankine--Hugoniot jump condition, see~\cite{DafermosBook}.
      \item Any discontinuity away from the constraint is classical, \textit{i.e.}~satisfies the Lax entropy inequalities, see~\cite{DafermosBook}.
      \item Nonclassical discontinuities, see~\cite{LeflochBook}, may occur only at the constraint location $x=0$, and in this case the flow at $x=0$ is the maximal flow allowed by the constraint. Namely, if the solution contains a nonclassical discontinuity for all times $t \in I$, $I$ open in $\reals_+$, then for a.e.~$t$ in $I$
          \begin{align}\label{eq:ws}
            f\left( \rho(t, 0-) \right) = f\left( \rho(t, 0+) \right) = p\left( \int_{\reals_-} w(x) ~ \rho(t,x) ~{\d} x\right) ~.
          \end{align}
    \end{enumerate}
\end{proposition}
If the constraint function $p$ is  multivalued the equalities in~\ref{eq:constrianed4} and~\ref{eq:ws} should be interpreted as inclusions, and the result of the proposition remains true.

The existence result in \cite{AndreianovDonadelloRosini} is achieved by a procedure which couples the operator splitting method \cite{DafermosHsiao},  with the wave--front tracking algorithm, \cite{DafermosWFT}, see also \cite{AMADORISHEN} for a similar technique. This approach allows us to approximate our problem by a problem with ``frozen'' constraint, as~\ref{eq:constrianed}, at each time step.

The regularity of $p$ plays a central role in the well-posedness result.  While existence still holds for the Cauchy problem when $p$ is merely continuous, it is difficult to justify uniqueness in this case. Further,  in \cite{AndreianovDonadelloRosini}, the authors give some basic examples illustrating that solutions of a Riemann problem for the case of a non-decreasing piecewise constant $p$, see \textbf{(P2)}  below, may fail to be unique, $\Lloc1$--continuous and consistent.

The case in which $p$ is piecewise constant is extremely important both for the theoretical study of the problem and its numerical applications. First, it is related to the construction of the Riemann solver, which is the basic building block for the wave--front tracking algorithm, a precious tool in the study of existence and stability of the solutions for the general Cauchy problem. Moreover, the piecewise constant case is essentially the only case in which solutions can be computed explicitly, which is an undeniable \emph{atout} when looking for examples and applications.  To this aim, it is relevant to provide a detailed study of the different pathological behaviors one may encounter. Remarkably, we show that these behaviors can be easily forecast and avoided when looking for explicit examples of solutions.

In the following section we develop a detailed proof of the fact that, when $p$ is piecewise constant, the Riemann solver for~\ref{eq:constrianed} with a constraint of the form~\ref{eq:constrianed4}, is not unique and does not satisfy the minimal requirements needed to develop the classical wave--front tracking approach. Additionally, we compare the two extreme Riemann solvers: the one that minimizes the capacity drop, and the one that maximizes it. In particular, for any time $T>0$ we estimate the distance between the profiles of the solutions produced by the two Riemann solvers starting from the same initial condition.

\section{The constrained Riemann problem}\label{sec:Riemann}

In this section we study constrained Riemann problems of the form
\begin{subequations}\label{eq:constrianedRiemann}
\begin{align}\label{eq:constrianedRiemann1}
    \partial_t\rho + \partial_xf(\rho) &= 0 & (t,x) &\in \reals_+\times\reals\\
    \label{eq:constrianedRiemann2}
    f\left(\rho(t,0\pm)\right) &\le p\left( \int_{\reals_-} w(x) ~ \rho(t,x) ~{\d} x \right) & t &\in \reals_+\\
    \label{eq:constrianedRiemann3}
    \rho(0,x) &=\left\{
    \begin{array}{l@{\qquad\hbox{if }}l}
      \rho_L & x < 0\\
      \rho_R & x \ge 0
    \end{array}
    \right. & x &\in \reals~,
\end{align}
\end{subequations}
 with $\rho_L, \rho_R \in [0,R]$. The flux $f$ and the weight function $w$ satisfy~\textbf{(F)} and~\textbf{(W)}, moreover,  we adopt the following assumption on $p$ (instead of~\textbf{(P1)}) to allow an explicit construction of solutions to~\ref{eq:constrianedRiemann}
\begin{enumerate}
  \item[\textbf{(P2)}] $p ~\colon~ \left[0,R\right] \to \left]0,f(\bar\rho)\right]$ is piecewise constant non--increasing map with a finite number of jumps.
\end{enumerate}
\begin{figure}[htpb]
\centering
        \includegraphics[width=.75\textwidth]{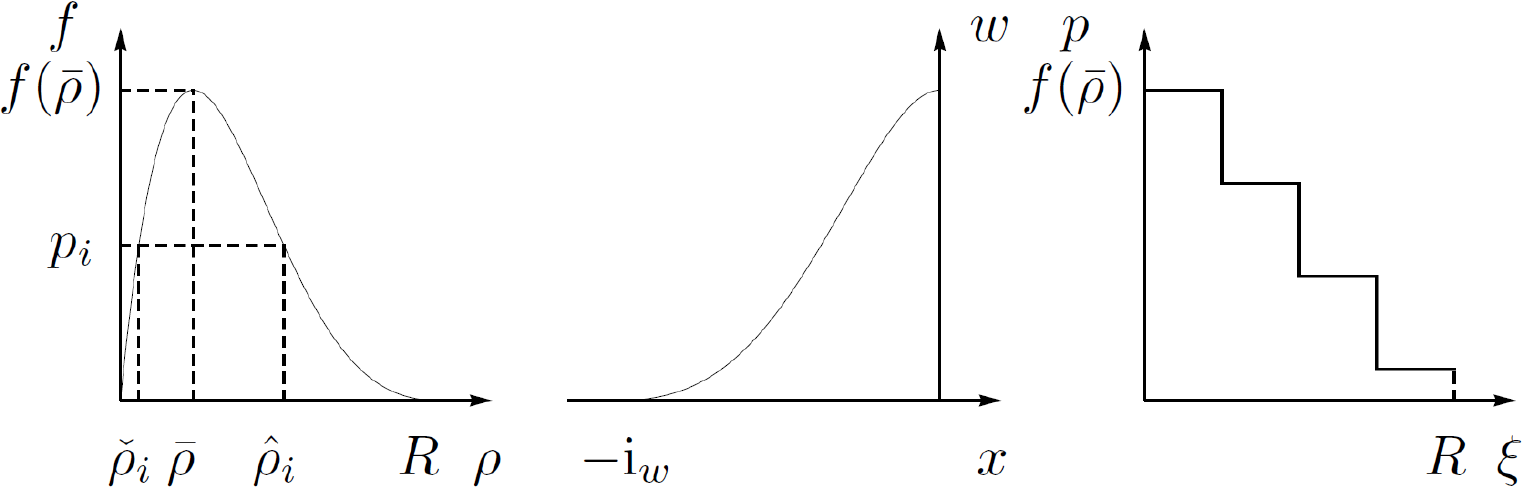}
      \caption{Examples of functions satisfying conditions~\textbf{(F)}, \textbf{(W)},  and \textbf{(P2)}.}
\label{fig:Korn}
\end{figure}
 Unfortunately, the regularity of $p$ required by~\textbf{(P2)} is not enough to apply the well--posedness results in \cite{AndreianovDonadelloRosini}. In particular, Example~2 in \cite{AndreianovDonadelloRosini} illustrates the loss of uniqueness and stability of entropy weak solutions. In this section we present a systematic study of the possible pathological behaviors. We denote by $\mathcal{R}$ the classical Riemann solver. This means that the map $\left[ (t,x) \mapsto \mathcal{R}[\rho_L, \rho_R](x/t) \right]$ is the unique entropy weak solution for the unconstrained problem~\ref{eq:constrianedRiemann1}, \ref{eq:constrianedRiemann3}, see for example~\cite{BressanBook} for its construction. Whenever the classical weak solution given by $\mathcal{R}$ does not satisfy the constraint~\ref{eq:constrianedRiemann2}, we replace it by a nonclassical weak solution,  see~\cite{LeflochBook} as a general reference,
\begin{equation}
    \rho(t,x) = \left\{
    \begin{array}{l@{\qquad\hbox{if }}l}
      \mathcal{R}[\rho_L,\hat\rho( q)](x/t)& x<0\\
      \mathcal{R}[\check\rho(q),\rho_R](x/t)& x\ge0 ~,
    \end{array}
    \right.
\end{equation}
where  the maps $\check\rho,\hat\rho ~\colon~ [0,f(\bar\rho)] \to [0,R]$  are implicitly defined by
\begin{align*}
    f\left(\check\rho(q)\right) = q = f\left(\hat\rho(q)\right)\quad\hbox{ and }\quad\check\rho(q) \le\bar\rho \le \hat\rho(q)~.
\end{align*}
\begin{figure}[htpb]
\centering
        \includegraphics[width=\textwidth]{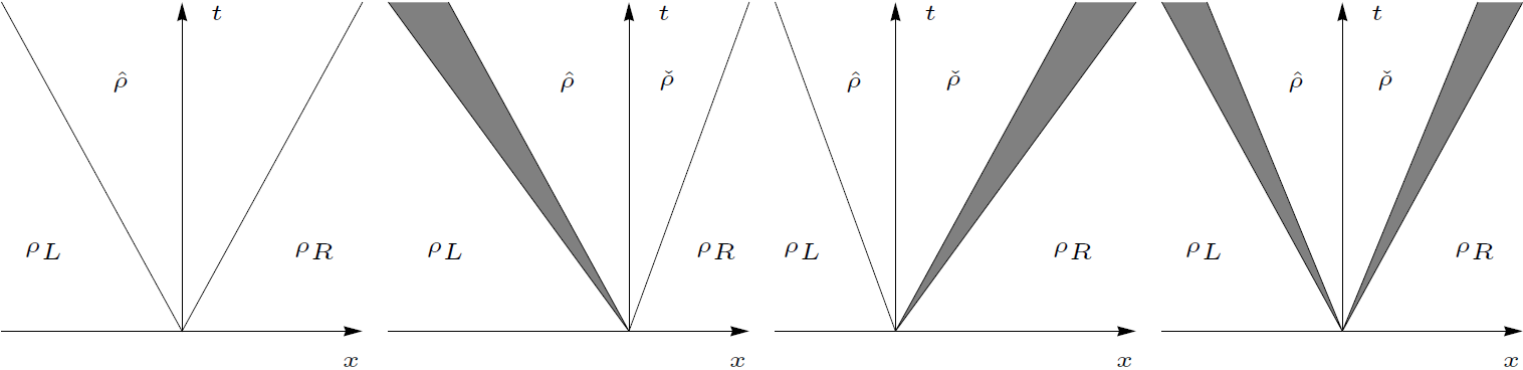}
      \caption{The four possible configurations of nonclassical entropy weak solutions of the form~\ref{eq:nonclassicalsol}.}
\label{fig:N123}
\end{figure}
We stress that by Proposition~\ref{prop:ws} any nonclassical entropy weak solution  is a classical entropy weak solution in the Kru\v zkov sense, see \cite{Kruzkov}, \cite{DafermosBook}, in the half--planes $\reals_+ \times \reals_-$ and $\reals_+ \times \reals_+$. However, the jump at $x= 0$ is a nonclassical shock, in the sense that it does not satisfy the Lax entropy inequalities.

First, we should notice that as soon as the constraint function $t \to q(t)$ is not constant, the solutions of the Riemann problem may not be self--similar.
\begin{example}
  We use a Cauchy problem of the form~\ref{eq:constrianed} to model vehicular traffic in presence of a traffic light. Assume $f(\rho) = \rho(1-\rho)$ and $q(t) = 0.25 \sum_{k\in\mathbb{N}} \chi_{[2k,2k+1[}(t)$. When the traffic light is green, \emph{i.e.}~for $t\in[2k,2k+1[$, the flow at $x=0$ is free from any constraint. Conversely, when the traffic light is red, \emph{i.e.}~for $t\in\left[2k+1,2(k+1)\right[$, the admissible flow at $x=0$ become zero. This means that the Riemann problem with initial condition at $t=0$ given by $\rho_L = \rho_R \neq 0$ will not be self--similar, because the constant solution will not satisfy the constraint starting from $t= 1$.

The above example also shows that as soon as we use a nonclassical Riemann solver we lose the \emph{a priori} $\BV$ bounds on the solution.
\end{example}

In the proof of Proposition~\ref{prop:riemann} we show that any entropy weak solution of~\ref{eq:constrianedRiemann} is self--similar for sufficiently small times. Therefore, it makes sense to introduce nonclassical \emph{local} Riemann solvers, see Definition~\ref{def:Rp}. Then, the availability of a local Riemann solver allows us to construct a \emph{global} solution to the Riemann problem~\ref{eq:constrianedRiemann} by a wave--front tracking algorithm in which the jumps in the map $[t \mapsto q(t)]$ are interpreted as non--local interactions.

Aiming for a general construction of the solutions to~\ref{eq:constrianedRiemann}, we allow $p$ to be a multi--valued piecewise constant function, namely, see Fig.~\ref{fig:Korn}, right:
\begin{itemize}
\item there exist $\xi_1, \ldots, \xi_n \in \left]0,R\right[$ and $p_0, \ldots, p_n \in \left]0,f(\bar\rho)\right]$, with $\xi_i < \xi_{i+1}$ and $p_i > p_{i+1}$, such that $p(0) = p_0$, $p(R) = p_n$, $p ~\caratt{\left]\xi_i, \xi_{i+1}\right[} = p_i$ for $i=0,\ldots,n-1$,  $p(\xi_i)= \left[p_i, p_{i-1}\right]$ for $i=1,\ldots,n$, being $\xi_0 = 0$ and $\xi_{n+1} = R$.
\end{itemize}
In the following we will use the notations $\check\rho_i = \check\rho(p_i)$ and $\hat\rho_i = \hat\rho(p_i)$.

As it will become clear in Proposition~\ref{prop:riemann},  the possible loss of uniqueness and stability can be easily forecast once the piecewise constant constraint $p$ and  the flux $f$ are given. In particular, for some respective configurations of $p$ and $f$ the solution of the Riemann problem exists and is unique, locally in time, for any initial data $\left(\rho_L, \rho_R\right)$ in $\left[0, R\right]^2 $.

\begin{definition}\label{def:regionsRP}
  Introduce the subset of $\left[0, R\right]^2$
  \begin{align*}
    \mathcal{C} &= \left\{ \left(\rho_L, \rho_R\right) \in \left[0, R\right]^2 ~\colon~ \left(\rho_L, \rho_R\right)\hbox{satisfies condition~\textbf{(C)}} \right\} ,
  \end{align*}
  where we say that $\left(\rho_L, \rho_R\right)$ satisfies condition~\textbf{(C)} if it satisfies one of the following conditions:
  \begin{description}
  \item[(C1)] $\rho_L < \rho_R$, $f(\rho_R) < f(\rho_L)$ and $f(\rho_R) \le p(\rho_L+)$;
  \item[(C2)] $\rho_L <\rho_R$, $f(\rho_L) \le f(\rho_R)$ and $f(\rho_L) \le p(\rho_L+)$;
  \item[(C3)] $\rho_R \le \rho_L \le \bar\rho$ and $f(\rho_L) \le p\left( \rho_L +\right)$;
  \item[(C4)] $\rho_R \le \bar\rho < \rho_L$ and $f(\bar\rho) = p\left( \rho_L +\right)$;
  \item[(C5)] $\bar\rho < \rho_R \le \rho_L$, $f(\rho_R) \le p\left( \rho_L -\right)$ and $f(\rho_L) < p\left( \rho_L +\right)$.
  \end{description}

  Analogously, introduce the subset of $\left[0, R\right]^2$
  \begin{align*}
    \mathcal{N}  &= \left\{ \left(\rho_L, \rho_R\right) \in \left[0, R\right]^2 ~\colon~ \left(\rho_L, \rho_R\right)\hbox{satisfies condition~\textbf{(N)}} \right\},
  \end{align*}
  where we say that $\left(\rho_L, \rho_R\right)$ satisfies condition~\textbf{(N)} if it satisfies one of the following conditions:
  \begin{description}
  \item[(N1)] $\rho_L < \rho_R$ and $f(\rho_L)> f(\rho_R) > p(\rho_L+)$;
  \item[(N2)] $\rho_L <\rho_R$, $f(\rho_L) \le f(\rho_R)$ and $f(\rho_L)> p(\rho_L-)$;
  \item[(N3)] $\rho_R \le \rho_L \le \bar\rho$ and $f(\rho_L) > p\left( \rho_L -\right)$;
  \item[(N4a)] $\rho_R \le \bar\rho < \rho_L$, $f(\bar\rho) \ne p\left( \rho_L- \right)$ and $f(\rho_L) < p(\rho_L+)$;
  \item[(N4b)] $\rho_R \le \bar\rho < \rho_L$, $f(\bar\rho) \ne p\left( \rho_L- \right)$ and $f(\rho_L)> p(\rho_L-)$;
  \item[(N5a)] $\bar\rho < \rho_R \le \rho_L$, $f(\rho_R) > p\left( \rho_L- \right)$ and $f(\rho_L) < p\left( \rho_L +\right)$;
  \item[(N5b)] $\bar\rho < \rho_R \le \rho_L$ and $f(\rho_L) > p\left( \rho_L -\right)$.
  \end{description}
\end{definition}

Here, $\mathcal C$ stands for \emph{classical} and $\mathcal N$, for \emph{nonclassical}, in relation with the nature of the shock appearing in the solution of~\ref{eq:constrianedRiemann} at $x=0$.
Observe that if the constraint function $p$ is constant in a neighborhood of the state $\rho_L$, then $p(\rho_L-) = p(\rho_L+)$ and this simplifies the above conditions. Also a right or left continuity assumption on $p$ would simplify the above definitions.

The next proposition says that uniqueness holds at least for small times if and only if the initial data are in $\mathcal{C} \cup \mathcal{N}$.  It is fundamental to remark that since non--uniqueness is possible only when  $p(\rho_L-) \ne p(\rho_L+)$ and $p(\rho_L-) \ge f(\rho_L) \ge p(\rho_L+)$, non--uniqueness concerns at most a finite number of left states and therefore the region $[0,R]^2 \setminus (\mathcal{C} \cup \mathcal{N})$ is the union of a finite number of line segments.

\begin{proposition}\label{prop:riemann}
    Consider the constrained Riemann problem~\ref{eq:constrianedRiemann}.

    \noindent$\bullet$~If $(\rho_L, \rho_R) \in \mathcal{C}$, then the map $\left[ (t,x) \mapsto \mathcal{R}[\rho_L, \rho_R](x/t) \right]$ is the unique entropy weak solution at least  for $t>0$ sufficiently small.

    \noindent$\bullet$~If $(\rho_L, \rho_R) \in \mathcal{N}$, then there exists a unique $\bar p \in \left[p(\rho_L+), p(\rho_L-)\right]$ such that the map
            \begin{align*}
                \left[
                t \mapsto \left\{\begin{array}{l@{\quad\hbox{if }}l} \mathcal{R}[\rho_L, \hat\rho\left( \bar p \right)](x/t)&x<0\\ \mathcal{R}[\check\rho\left( \bar p \right), \rho_R](x/t)&x\ge0\end{array}\right.\right]
            \end{align*}
            is the unique entropy weak solution at least for $t>0$ sufficiently small.

    \noindent$\bullet$~If $(\rho_L, \rho_R) \in [0,R]^2 \setminus (\mathcal{C} \cup \mathcal{N})$, then the corresponding constrained Riemann problem~\ref{eq:constrianedRiemann} admits more than one entropy weak solution.
\end{proposition}
The proof is deferred to Section~\ref{sec:proofRiem}

\begin{remark}
Once the function $p$ is fixed, the time interval $[0,\tau]$ on which the solution to the Riemann problem~\ref{eq:constrianedRiemann} is self-similar can be estimated from the initial datum $(\rho_L,\rho_R)$ provided it belongs to $\mathcal C \cup \mathcal N$.

On the contrary, when $(\rho_L, \rho_R) \in [0,R]^2 \setminus (\mathcal{C} \cup \mathcal{N})$, we are not always able to forecast when the next ``interaction with the constraint'' will take place. In some situations, a whole one--parameter family of solutions exists, we refer to Example 2 in ~\cite{AndreianovDonadelloRosini} for a detailed discussion of this case.
\end{remark}

\subsection{Riemann solvers}\label{sec:solvers}

As the local solutions of the Riemann problem~\ref{eq:constrianedRiemann} are not unique in general, we are naturally led to question the existence of suitable selection criteria. All the solutions we introduce are solutions in the Kru\v zkov sense in the open half--planes $\reals_+\times \reals_+$  and  $\reals_+\times \reals_-$, so they satisfy the basic requirement of entropy dissipation. However, coming back to the real situations which our model aims to describe, we argue that the most interesting behaviors to track correspond to the extreme cases in which the flux at the exit is either the highest or the lowest possible from a given initial condition.

If, as an example, our model describes the evacuation of a narrow corridor, it is clear that the optimal solution corresponds to the highest admissible values of the flux at the exit. By opposition to the next case, we describe this situation as \emph{quiet} behaviour. In analogy to the discussion in~\cite{Libbano} we interpret all other possible solutions as consequences of an irrational behavior, which in literature is often described as \emph{panic}. In particular, we can use the solution corresponding to the lowest admissible values of the flux at the exit to find an upper bound for the evacuation time.

From now on we restrict ourselves to the case in which $p(\xi_i)$ can only take the values $p_i$ and $p_{i+1}$ and not the intermediate values.

\begin{definition}\label{def:Rp}
    Two Riemann solvers $\mathcal{R}^q$ and $\mathcal{R}^p$ for~\ref{eq:constrianedRiemann} are defined as follows for $t>0$ sufficiently small and $x \in \reals$:
    \begin{description}
    \item[(C)] If $(\rho_L, \rho_R) \in \mathcal{C}$ then
      \[\mathcal{R}^q[\rho_L, \rho_R] (t,x)=\mathcal{R}^p[\rho_L, \rho_R](t,x) = \mathcal{R}[\rho_L, \rho_R](x/t)~.\]
    \item[(N)] If $(\rho_L, \rho_R) \in \mathcal{N}$ then
        $$\mathcal{R}^q[\rho_L, \rho_R] (t,x)= \mathcal{R}^p[\rho_L, \rho_R](t,x) = \left\{\begin{array}{l@{\quad\hbox{if }}l} \mathcal{R}[\rho_L, \hat\rho\left(\bar p \right)](x/t)&x<0\\ \mathcal{R}[\check\rho\left( \bar p \right), \rho_R](x/t)&x\ge0~,\end{array}\right.$$
      where $\bar p = p(\rho_L-)$ if $(\rho_L,\rho_R)$ satisfies~\textbf{(N4a)} or~\textbf{(N5a)}, otherwise $\bar p = p(\rho_L+)$.
    \item[(CN2), (CN3), (CNN5)] If $(\rho_L, \rho_R)$ satisfies one of these sets of conditions  then
        \begin{align*}
            \mathcal{R}^q[\rho_L, \rho_R](t,x) &= \mathcal{R}[\rho_L, \rho_R](x/t)~,\\
            \mathcal{R}^p[\rho_L, \rho_R](t,x) &= \left\{\begin{array}{l@{\quad\hbox{if }}l} \mathcal{R}[\rho_L, \hat\rho\left(p(\rho_L+)\right)](x/t)&x<0\\ \mathcal{R}[\check\rho\left(p(\rho_L+) \right), \rho_R](x/t)&x\ge0~.\end{array}\right.
        \end{align*}
     \item[(NNN4), (NNN5)] If $(\rho_L, \rho_R)$ satisfies one of these sets of conditions then  $\mathcal{R}^q[\rho_L, \rho_R](t,x) $ takes the form~\ref{eq:nonclassicalsol} with $\bar p = p(\rho_L-) $ and $\mathcal{R}^p[\rho_L, \rho_R](t,x) $ takes the form~\ref{eq:nonclassicalsol} with $\bar p = p(\rho_L+) $.
    \end{description}
    \end{definition}
In the next proposition we collect the main properties of the Riemann solvers $\mathcal{R}^q$ and $\mathcal{R}^p$. In particular \textbf{(R6)} means that the Riemann solver $\mathcal{R}^q$ is the one which allows for the fastest evacuation, while $\mathcal{R}^p$ is associated to the slowest one.
\begin{proposition}\label{prop:Riemann}
Let $(\rho_L, \rho_R) \in \left[0,R\right]^2$. Then, for $\star = q,\,p$:
\begin{enumerate}
\item[\textbf{(R1)}] $\left[(t,x) \mapsto \mathcal{R}^\star[\rho_L, \rho_R](t,x) \right]$ is a weak solution to~\ref{eq:constrianedRiemann1}, \ref{eq:constrianedRiemann3}.
\item[\textbf{(R2)}] $\mathcal{R}^\star[\rho_L, \rho_R]$ satisfies the constraint~\ref{eq:constrianedRiemann2} in the sense that
    \begin{align*}
        f\left(\mathcal{R}^\star[\rho_L, \rho_R](t,0\pm)\right) &\le p\left( \int_{\reals_-} w(x) ~ \mathcal{R}^\star[\rho_L, \rho_R] \left(t,x\right) ~{\d} x\right) .
    \end{align*}
\item[\textbf{(R3)}] $\mathcal{R}^\star[\rho_L, \rho_R](t) \in \BV\left( \reals;[0,R] \right)$.
\item[\textbf{(R4)}] The map $\mathcal{R}^\star ~\colon~ [0,R]^2 \to \Lloc1(\reals_+ \times \reals; [0,R])$ is continuous in $\mathcal{C}\cup\mathcal{N}$ but not in all $[0,R]^2$.
\item[\textbf{(R5)}] $\mathcal{R}^\star$ is consistent, see~\cite{ColomboGoatinConstraint}, \cite{ColomboPriuli} and the comment below.
\item[\textbf{(R6)}] $\mathcal{R}^q[\rho_L,\rho_R]$ maximizes the flux at the exit, in the sense that if $\mathcal{E}$ is the set of all entropy weak solutions of the Riemann problem~\ref{eq:constrianedRiemann}, we have
  $$
  \max_{\rho \in \mathcal{E}}\left\{f(\rho(t,0\pm)) \right\} = f\left(\mathcal{R}^q[\rho_L,\rho_R](0\pm)\right).
  $$
  Analogously, $\mathcal{R}^p[\rho_L,\rho_R]$ minimizes the flux at the exit, in the sense that
  $$
  \min_{\rho \in \mathcal{E}}\left\{f(\rho(t,0\pm)) \right\} = f\left(\mathcal{R}^p[\rho_L,\rho_R](0\pm)\right).
  $$
\end{enumerate}
\end{proposition}

We recall that a Riemann solver is said to be consistent when the juxtaposition of the solutions of two Riemann problems with respective initial conditions $(\rho_L,\rho_M)$ and $(\rho_M,\rho_R)$ is the solution of the Riemann problem with datum  $(\rho_L,\rho_R)$. Moreover the vice versa also holds true, in the sense that whenever the state $\rho_M$ is an intermediate state in the solution of the Riemann problem with initial condition  $(\rho_L,\rho_R)$, then the solution consist of exactly the same states and waves which we would obtain by solving side by side the two Riemann problems with data $(\rho_L,\rho_M)$ and $(\rho_M,\rho_R)$.

 The proof of Proposition~\ref{prop:Riemann} is deferred to Section~\ref{sec:technicalRiemann}.

 \begin{remark}
   It is important to observe that even if $p(\xi_i)$ can only take the two values $p_i$ and $p_{i+1}$, this is not enough to rule out the existence of  infinitely many different solutions as the ones described in Example 2 of~\cite{AndreianovDonadelloRosini}, in the case $p_i > f(\xi_i) = p_{i+1}$. However, each of the extremes Riemann solvers $\mathcal{R}^\star$, $\star=p, \,q$, selects one of them because it sticks to the constant level of constraint prescribed by Definition~\ref{def:Rp}, the level $p_i$ for $\mathcal{R}^q$ and the level $p_{i+1}$ for $\mathcal{R}^p$, until a non--local interaction takes place.
 \end{remark}

\subsection{On the use of Riemann solvers \texorpdfstring{$\mathcal{R}^\star$, $\star = q,\,p$}{\emph{quiet} and \emph{panic}}}

Although the Riemann solvers $\mathcal{R}^\star$ are not $\Lloc1$--continuous, an existence result for the Cauchy problem~\ref{eq:constrianed} can be obtained from a wave--front tracking algorithm based on $\mathcal{R}^\star$, see for instance~\cite{ColomboRosini2}, \cite{Rosini09}. Such approach using $\mathcal{R}^\star$ does not require the operator splitting method. However, the non--local nature of the approximating problems prevents us from a direct application of the Riemann solvers $\mathcal{R}^\star$. In fact, even in a arbitrary small neighborhood of $x=0$, to prolong the approximating solution $\rho^n$ beyond a time $t = \bar t>0$ it is not sufficient to know the traces $\rho^n(\bar t,0-)$, $\rho^n(\bar t,0+)$, but also the value $\int_{-\iw}^0 w(x) ~\rho^n(\bar t,x) ~{\d}x$ is needed. Roughly speaking, because of the non--local character of the constraint one cannot merely juxtapose the solution to the Riemann problem associated to the values of the traces at $x=0$ with the solution to the Riemann problems away from the constraint. Finally, also jumps in $[t \mapsto p\left(\xi(t)\right)]$ have to be considered as (non--local) interactions. Therefore, the approach using $\mathcal{R}^\star$ is considerably heavier and more technical than the one presented in \cite{AndreianovDonadelloRosini}.

\subsection{On the comparison between the Riemann solvers \texorpdfstring{$\mathcal{R}^p$ and $\mathcal{R}^q$}{\emph{quiet} and \emph{panic}}}\label{subsec:twosolvers}

In this section we aim to compare the solutions obtained by the two Riemann solvers introduced above, starting from the same initial condition  $(\rho_L, \rho_R)$. It is clear from Definition~\ref{def:Rp} that as soon as   $(\rho_L, \rho_R) $ belongs to $ \mathcal{C}\cup\mathcal{N}$ the solutions obtained by the two Riemann solvers coincide.

As a preliminary remark we stress that adapting the proof of Theorem 3.1 of~\cite{AndreianovDonadelloRosini} to the case in which $p$ is discontinuous,  one can get a rough upper bound, exponential in time, for the $\L1$ distance between two solutions $\rho_1$, $\rho_2$ obtained from the same initial condition (not necessarily of Riemann type). Indeed, instead of the bound
\begin{equation}
    \modulo{p(\rho_1)-p(\rho_2)} \leq \lip(p) ~\modulo{\rho_1-\rho_2}~,
\end{equation}
valid when $p$ is Lipschitz continuous, in the discontinuous case one can use the bound
\begin{equation}\label{eq:rough-estimate-rho12}
    \modulo{p(\rho_1)-p(\rho_2)} \leq h + N ~\modulo{\rho_1-\rho_2}~,
\end{equation}
where $h$ is the maximal size of jump in $p(\cdot)$ and $N$ is a constant (observe that
if $p(\cdot)$ is a discretization of some Lipschitz function, see Section 4.1 of~\cite{AndreianovDonadelloRosini},
then $N$ can be taken independent of $h$).
From the fundamental stability estimate of Proposition~2.10 of~\cite{scontrainte},
using~\ref{eq:rough-estimate-rho12} and the Gronwall inequality one easily gets the bound
\begin{equation}\label{eq:rough-exponential-bound}
\norma{\rho_1(t)-\rho_2(t)}_{\L1(\reals;\reals)} \leq \frac{h}{2N} \left[\exp(2Nt)-1\right],
\end{equation}
whenever $\rho_1(0,\cdot)=\rho_2(0,\cdot)$. This rough estimate is enough to show that as $h$ goes to zero, the discrepancy between  different solutions vanishes and this argument applies to any initial datum, not necessarily of the Riemann type.

However, the exponential growth  with respect to $t$ of the  upper bound~\ref{eq:rough-exponential-bound}
is clearly not optimal when we aim to compare the solutions of a Riemann problem,  since it does not take into account the specific self-similar structure
of solutions valid at least on a small interval of time $[0,T]$.

Let us demonstrate that in the case where different Riemann solvers  co-exist,
the $\L1$ distance of the associated solution grows at most linearly both in $h$
and in $t\in [0,T]$ (see also the numerical experiment on Figure~\ref{rho_diff}).

\begin{figure}[htpb]
\centering
  \includegraphics[width=.31\textwidth]{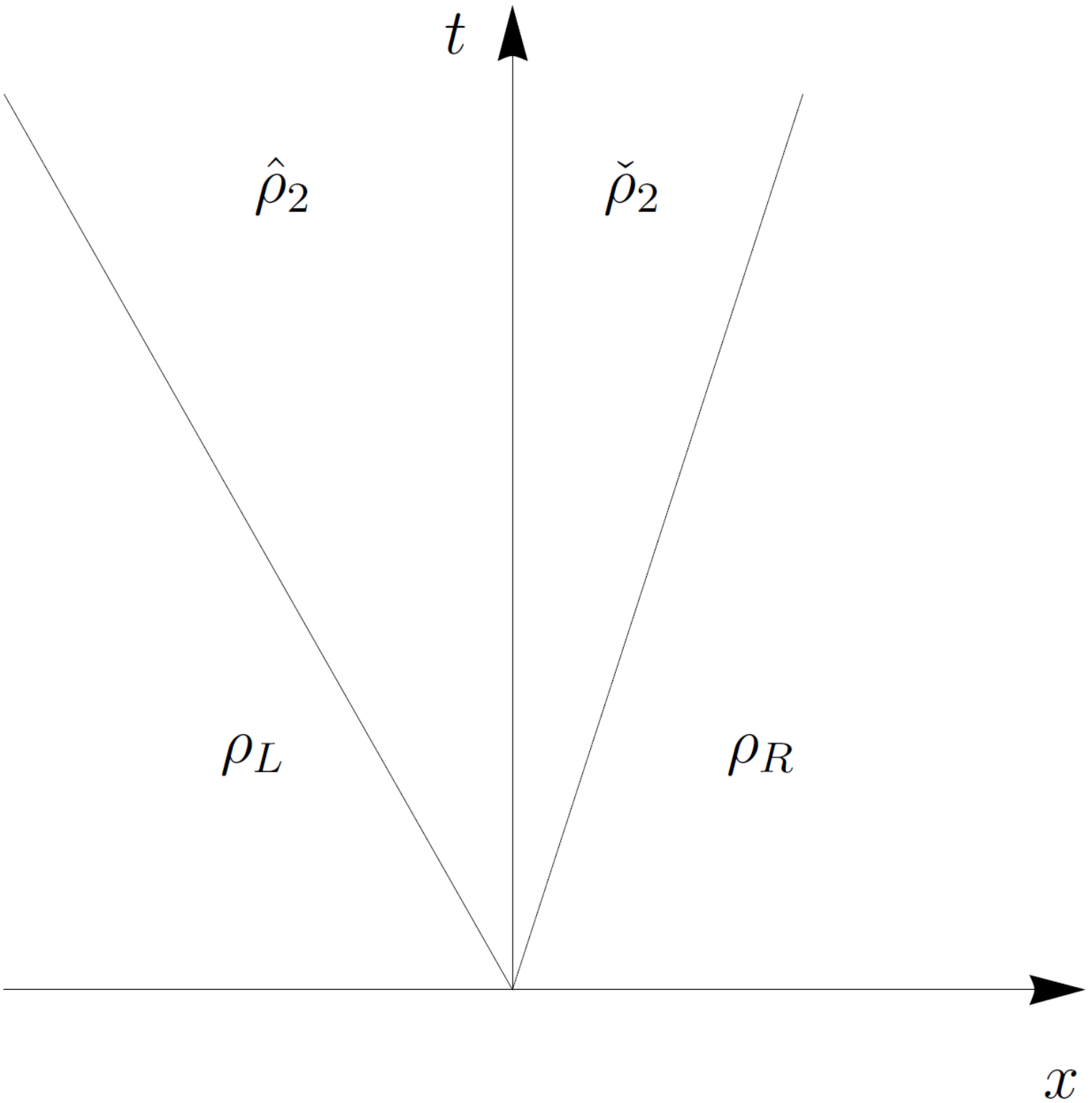}
  \includegraphics[width=.31\textwidth]{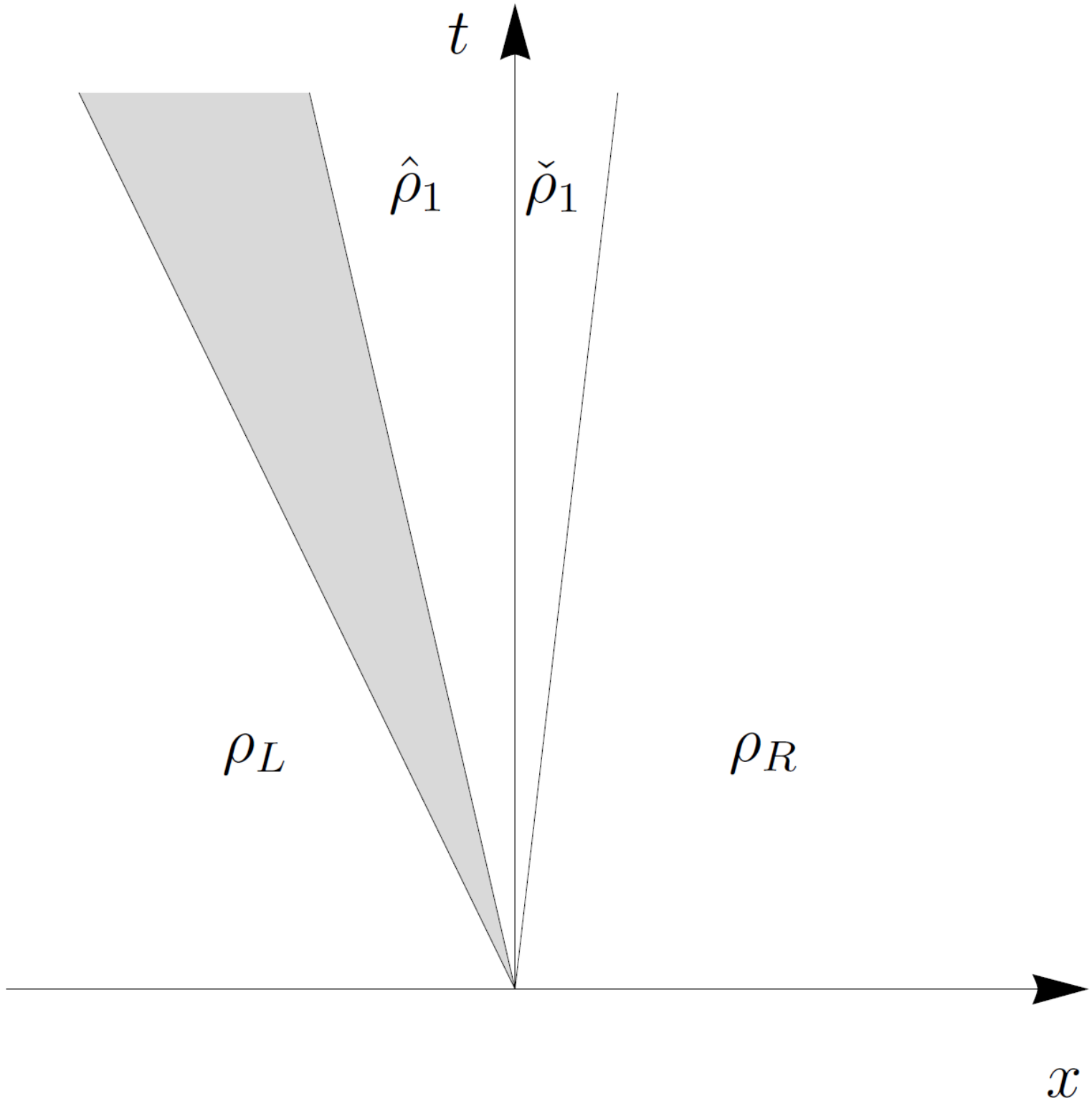}
  \includegraphics[width=.31\textwidth]{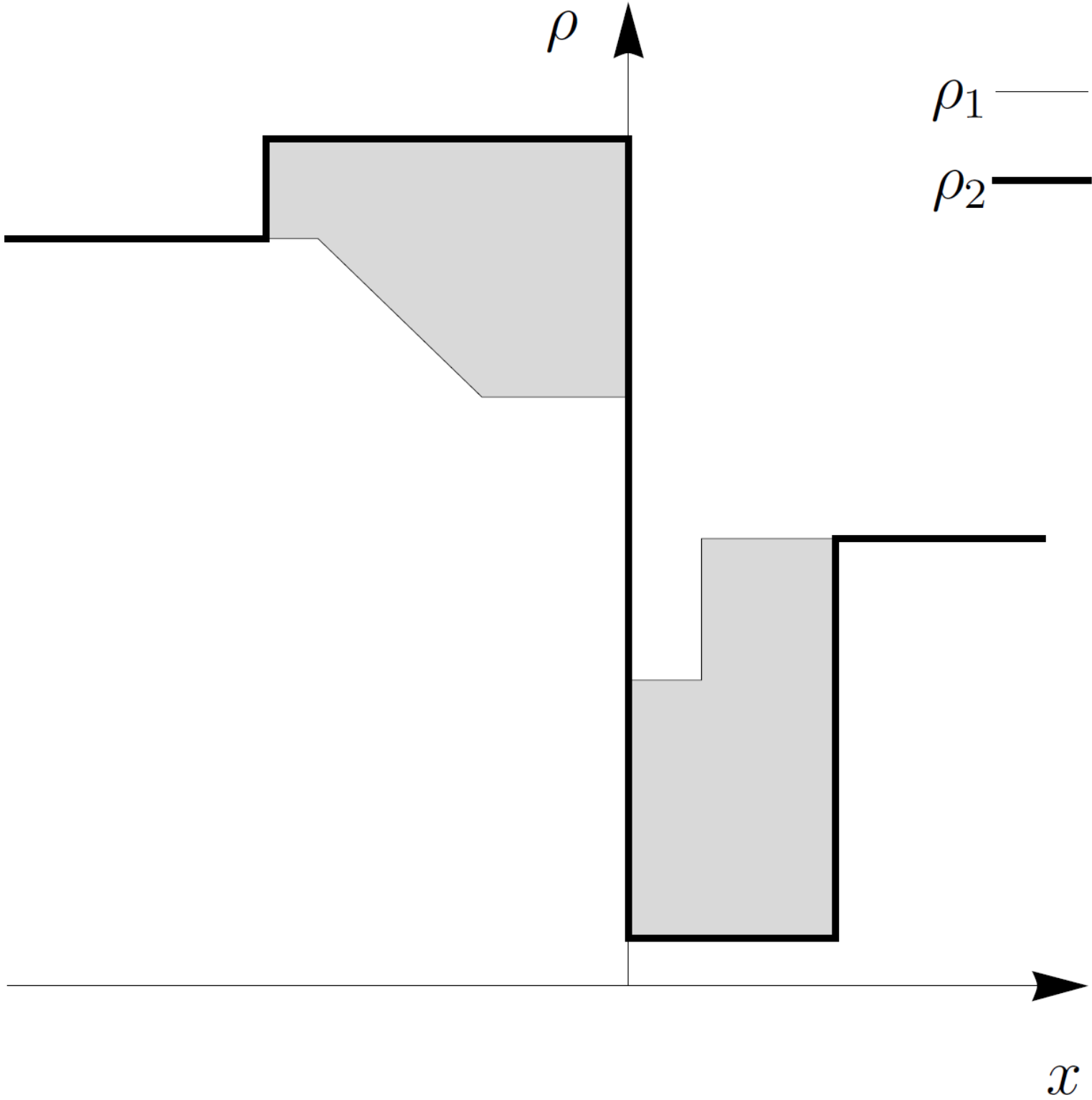}
  \caption{The solutions corresponding to $\mathcal{R}^p$, left, and to $\mathcal{R}^q$, center, and the comparison between their profiles at fixed time, right, as described in Section~\ref{subsec:twosolvers}.}
\label{fig:super}
\end{figure}

In order to keep our presentation as light as possible, we focus on only one of the possible cases in which $\mathcal{R}^p$ and $\mathcal{R}^q$ differ. All other cases can be handed in a similar way.
\begin{figure}[htpb]
\centering
  \includegraphics[width=.31\textwidth]{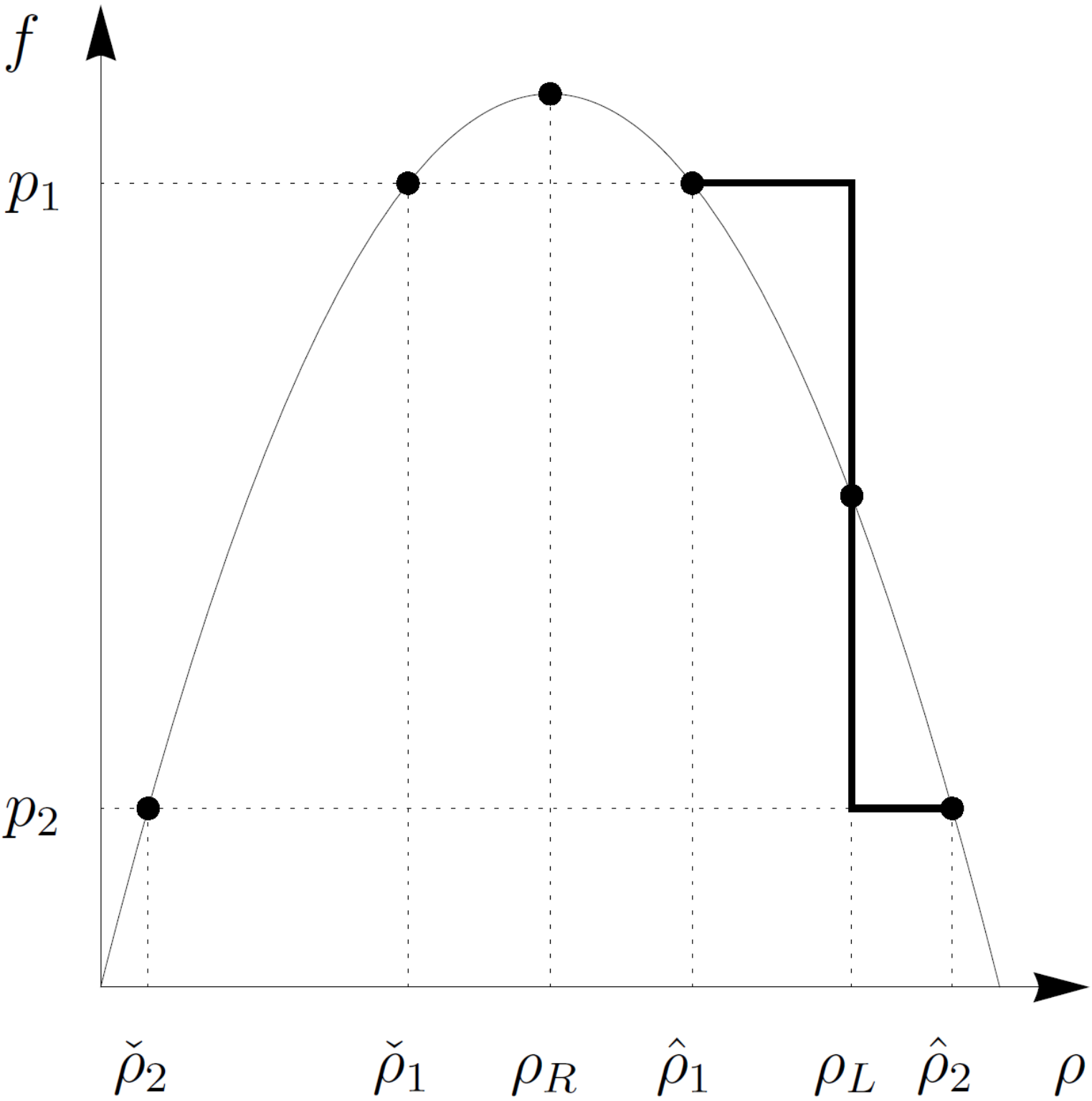}
  \caption{The flux $f$, the constraint function $p|_{]\hat\rho_1,\hat\rho_2[}$ and the values of the density $\rho$ considered in Section~\ref{subsec:twosolvers}.}
\label{fig:superflux}
\end{figure}
Assume that $p(\xi) = p_1 ~\chi_{[0,\bar\xi]}(\xi) + p_2 ~\chi_{]\bar\xi, R]}(\xi)$, where $\bar f>p_1>f(\bar\xi)>p_2>0$, and $\rho_L= \bar\xi$, $\rho_R = \bar\rho$, see Figure~\ref{fig:superflux}. We only consider solutions in a small interval of time $[0,T]$ in which they are self--similar. Then we get
\begin{align*}
    \mathcal{R}^p[\rho_L,\rho_R](t,x) &= \left\{
    \begin{array}{l@{\quad\hbox{if }}l}
    \mathcal{R}[\rho_L, \hat\rho_{2}](x/t) & x <0\\
    \mathcal{R}[\check\rho_{2}, \bar\rho](x/t) & x \ge0~,
    \end{array}
    \right.
    \\
    \mathcal{R}^q[\rho_L,\rho_R](t,x)&= \left\{
    \begin{array}{l@{\quad\hbox{if }}l}
    \mathcal{R}[\rho_L, \hat\rho_1](x/t) & x <0\\
    \mathcal{R}[\check\rho_1, \bar\rho](x/t) & x \ge0~,
    \end{array}
    \right.
\end{align*}
where the values $\hat\rho_{i}$ and $\hat\rho_{i}$, for $i=1,\,2$, are implicitly defined by the relations $f(\hat\rho_{i}) = f(\check\rho_{i}) = p_i$ and $\check\rho_{i}\leq \bar\rho\leq\hat\rho_{i}$. More explicitly, we can say that the solution corresponding to $\mathcal{R}^p$ consists of a shock of negative speed $\lambda_p$ between $\rho_L$ and $\hat\rho_{2}$, a stationary nonclassical shock between $\hat\rho_{2}$ and $\check\rho_{2}$ and a shock of positive speed $\mu_p$ between $\check\rho_{2}$ and $\bar\rho$. The solution corresponding to $\mathcal{R}^q$ consists of a rarefaction wave between $\rho_L$ and $\hat\rho_{1}$, a stationary nonclassical shock between $\hat\rho_{1}$ and $\check\rho_{1}$ and a shock of positive speed $\mu_q$ between $\check\rho_{1}$ and $\bar\rho$, see Figure~\ref{fig:super}.

As the characteristics of this problem propagate with finite speed, we expect the solutions associated to the two solvers coincide outside a bounded interval. The geometry of the problem, see Figure~\ref{fig:superflux}, implies that $\mu_p>\mu_q$ and that $\lambda_p$ is smaller than all the propagation speeds in the rarefaction wave between $\rho_L$ and $\hat\rho_{1}$. Therefore, at time $t\in[0,T]$ fixed, the two solutions coincide outside the interval $[\lambda_pt, \mu_pt]$.

The value of the distance $\norma{\mathcal{R}^p[\rho_L,\rho_R](t)-\mathcal{R}^q[\rho_L,\rho_R](t)}_{\L1(\reals;\reals)}$ corresponds, loosely speaking, to the value of the area between the profiles of solutions, see Figure~\ref{fig:super}.

Following the same technique as in \cite{BressanBook}, Chapter 7, we can estimate the distance between the solutions profiles at a fixed time $t\in[0,T]$. For the reader convenience, we recall that the propagation speed of the shock discontinuity between states the $\rho_a$ and $\rho_b$ is given by the Rankine-Hugoniot condition
\begin{equation}
  \sigma(\rho_a,\rho_b) = \frac{f(\rho_a) -f(\rho_b)}{\rho_a -\rho_b} ~,
\end{equation}
and that the propagation speed of the characteristics in a rarefaction wave joining the states $\rho_a$ and $\rho_b$ varies between the values $f'(\rho_a)$ and $f'(\rho_b)$. Also, by definition  $f(\hat\rho_i) = f(\check\rho_i) = p_i$, for $i=1,2$. A direct calculation gives us
\begin{align*}
    &\norma{\mathcal{R}^p[\rho_L,\rho_R](t,\cdot)-\mathcal{R}^q[\rho_L,\rho_R](t,\cdot)}_{\L1(\reals;\reals)}\\
    \leq& \norma{\mathcal{R}[\rho_L,\hat\rho_2](\cdot/t)-\mathcal{R}[\rho_L,\hat\rho_1](\cdot/t)}_{\L1(\reals;\reals)} +\norma{\mathcal{R}[\check\rho_2,\rho_R](\cdot/t)-\mathcal{R}[\check\rho_1,\rho_R](\cdot/t)}_{\L1(\reals;\reals)}\\
    \leq& \left[\lambda_p - \!f'(\rho_L)\right] \left(\hat\rho_2 - \rho_L\right)t + \!f'(\rho_L)\left(\hat\rho_2 - \hat\rho_1\right)t
    + \left(\check\rho_1-\check\rho_2\right)\mu_q \,t+ \left(\bar\rho - \check\rho_2\right)\left(\mu_q-\mu_p\right)t\\
    \leq& \,2~t ~\left[p_1-p_2 + o(\hat\rho_2-\hat\rho_1)\right].
\end{align*}
This means that whenever the piecewise constant function $p$ we consider is the discretization  of a smooth function we can bound \emph{a priori} the size of the error due to the lack of uniqueness and we can make it smaller and smaller as $h$ tends to $0$.

\section{Numerical results}

\begin{figure}[h!]
\centering
\includegraphics[width=0.8\hsize]{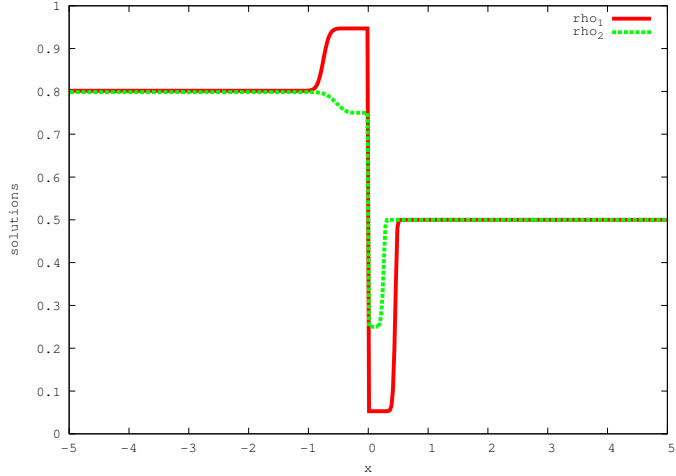}
\caption{The computed densities $\rho_1$ and $\rho_2$}
\label{Rp_Rq}
\end{figure}

\noindent We present here some numerical experiments in order to illustrate the results of the above section.
The scheme used for the simulations combines the ideas of~\cite{scontrainte} with an explicitly updated constraint computed
from weighted space averages of the discrete solution at previous time step. We will justify in the future work~\cite{ADRR}
convergence of this scheme to an entropy solution of the nonlocally constrainted problem in the sense of Definition~\ref{def:entropysol}, where the constraint function $p(\cdot)$ must be taken multi-valued. While it is delicate or even impossible to identify a unique Riemann solver to which the scheme would converge, we can use the simulations on Figure~\ref{Rp_Rq} to illustrate the fact that non-uniqueness for the Riemann problem results as unstable behavior in a vicinity of some specific data.

\noindent For the examples, we consider the flux $f(\rho)=\rho \left(1-\rho\right)$. The domain of computation is $x \in [-5,5]$, the constraint function is $p(\xi)=p_1 ~\chi_{[0,0.8]}(\xi)+p_2 ~\chi_{]0.8,1]}(\xi)$, where
$p_1=0.1875$, $p_2=0.05$, the weight function is $w(x)=2(x+1) ~\chi_{\left]-1,0\right]}(x)$. The final time of computation is $T=1$.
In Figure~\ref{Rp_Rq} is shown the computed solutions $\rho_1$ and $\rho_2$ corresponding respectively to the following initial states
\begin{align*}
    &\rho_0^1(x)=\left\{\begin{array}{l@{\quad\hbox{if }}l}
    0.8015&x<0\\
    0.5&x>0
    \end{array}\right.&\mbox{ and }&
    &\rho_0^2(x)=\left\{\begin{array}{l@{\quad\hbox{if }}l}
    0.7984&x<0\\
    0.5&x>0~.
    \end{array}\right.
\end{align*}
 \noindent Note that we have considered $\rho_L^1$ and $\rho_L^2$ such that  $\rho_L^1-\rho_L^2\simeq\Delta x/8$, where the space step $\Delta x=0.025$.
 Finally we assume the time step $\Delta t=\Delta x/10$.

\begin{figure}[h!]
\centering
\includegraphics[width=0.8\hsize]{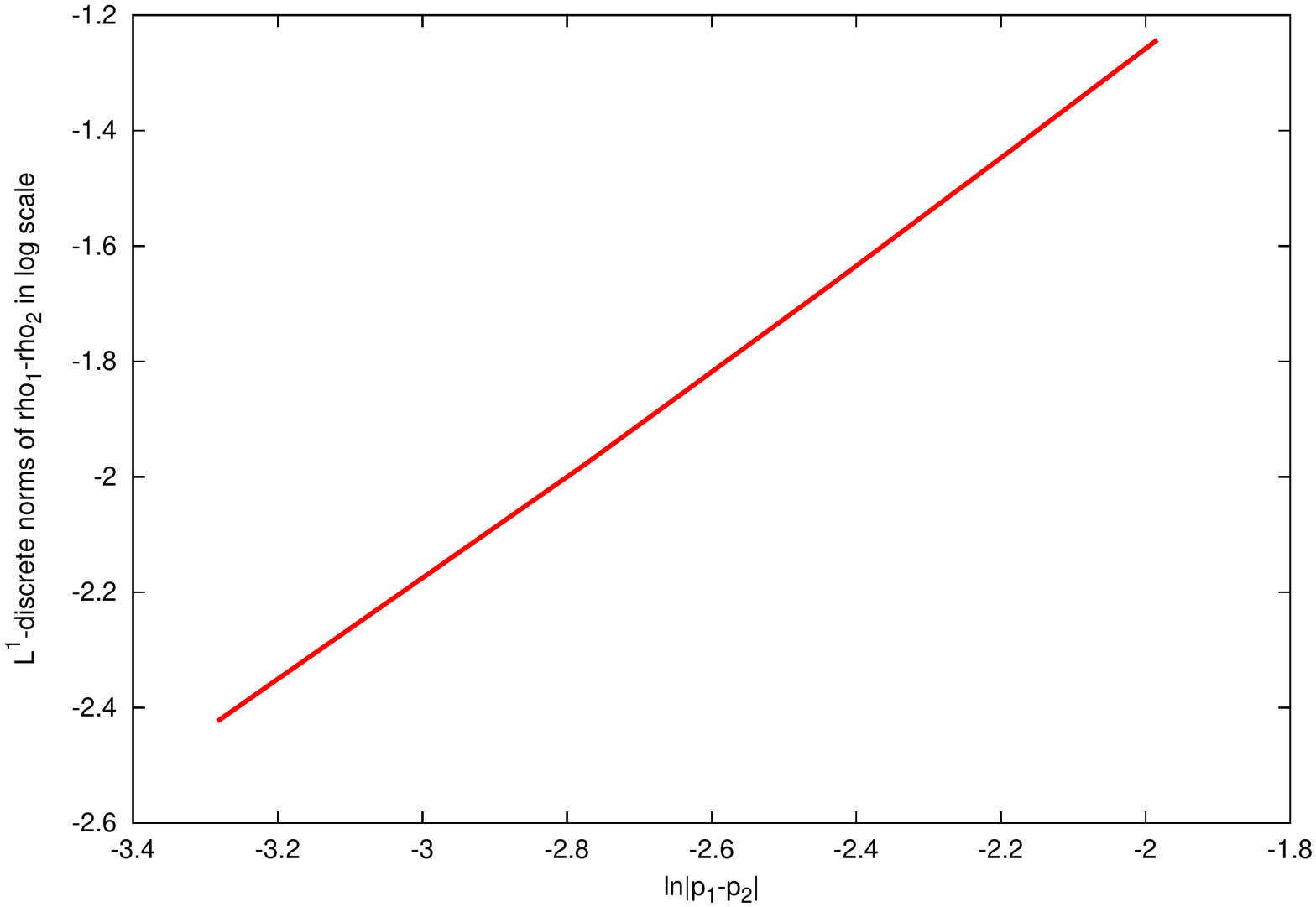}
\tiny\caption{The norm $\norma{\rho_1-\rho_2}_{\L1}$ with respect to $\modulo{p_1-p_2}$ in log/log scale.}
\label{rho_diff}
\end{figure}

\noindent
 Nonetheless, as shown in Figure~\ref{rho_diff}, in practice the instability is limited to a behavior
 of kind $$\norma{\rho_1-\rho_2}_{\L1([-5,5];\reals)} \le \norma{\rho_0^1-\rho_0^2}_{\L1([-5,5];\reals)}+C ~ h^\alpha ~,$$
 where $h$ is the maximal size of jump in $p(\cdot)$, $\alpha$ is close to $1$ and $C>0$ is a constant.

\noindent Indeed we reported in Figure~\ref{rho_diff} the computation of the $\L1$-discrete norms of the difference $\rho_1-\rho_2$ when we take, $p_2=0.05$, $0.075$, $0.1$, $0.125$ and $0.15$ in the definition of the constraint function. Using logarithmic scales, we deduce that the distance between the two solutions is approximatively proportional  to $\modulo{p_1-p_2}^{0.9}$.

\section{Proofs}

\subsection{Proof of Proposition~\ref{prop:riemann}}\label{sec:proofRiem}
First, we introduce the notation
\begin{align*}
    \xi(t) = \int_{\reals_-} w(x) ~\rho(t,x) ~{\d}x ~.
\end{align*}
Therefore  $\xi(0)=\rho_L$ and the map $[t \mapsto \xi(t)]$ is continuous.

We stress that any nonclassical entropy weak solution in the sense of Definition~\ref{def:entropysol} is also a classical entropy weak solution in the Kru\v zkov sense in the half--planes $\reals_+ \times \reals_-$ and $\reals_+ \times \reals_+$. Therefore, at least for $t>0$ sufficiently small, by Proposition~\ref{prop:ws}, assumption~\textbf{(P2)} and the continuity of the map $[t \mapsto \xi(t)]$, any nonclassical entropy weak solution of~\ref{eq:constrianedRiemann} must have the form, see Fig.~\ref{fig:N123},
\begin{subequations}\label{eq:nonclassicalsol}
\begin{align}\label{eq:nonclassicalsol1}
    \rho(t,x) = \left\{
    \begin{array}{l@{\qquad\hbox{if }}l}
      \mathcal{R}[\rho_L,\hat\rho(\bar p)](x/t)& x<0\\
      \mathcal{R}[\check\rho(\bar p),\rho_R](x/t)& x\ge0 ~.
    \end{array}
    \right.
\end{align}
Observe that~\ref{eq:nonclassicalsol1} is uniquely identified once we know $\bar p$ which, by~\ref{eq:ws}, satisfies
   \begin{align}\label{eq:nonclassicalsol2}
        \bar p = f\left(\check\rho(\bar p)\right) = f\left(\hat\rho(\bar p)\right) ~.
    \end{align}
We recall that~\ref{eq:nonclassicalsol2} means in particular that the Rankine--Hugoniot jump condition is satisfied at $x=0$ even when the solution to the Riemann problem is nonclassical. As a consequence of~\ref{eq:nonclassicalsol2}, of assumption~\textbf{(P2)} and of the continuity of $\left[t \mapsto \xi(t)\right]$, we have that
\begin{align}\label{eq:nonclassicalsol3}
    \bar p \in \left[p(\rho_L+), p(\rho_L-)\right].
\end{align}
\end{subequations}
This implies that if $p(\rho_L+) = p(\rho_L-)$, then $p(\xi)$ is constant in a neighborhood of $\rho_L$ and, since the solution is in $\C0\left(\reals_+; \Lloc1(\reals;[0,R])\right)$, uniqueness is ensured by the results in~\cite{ColomboGoatinConstraint}. However, the continuity of $p$ at $\rho_L$ is not a necessary condition for uniqueness as we show in the following section.

\subsubsection{Cases in which uniqueness holds}

In this section we prove that:
\begin{description}
  \item[If $(\rho_L, \rho_R) \in \mathcal{C}$]  the corresponding classical solution satisfies~\ref{eq:constrianedRiemann} for all $t>0$ sufficiently small and it is not possible to construct a different solution.
  \item[If $(\rho_L, \rho_R) \in \mathcal{N}$] the corresponding classical solution does not satisfy~\ref{eq:constrianedRiemann2}, and there exists a unique nonclassical solution that satisfies~\ref{eq:constrianedRiemann}.
\end{description}

We list here two basic properties which will be of great help in the following case by case analysis.

By assumption~\textbf{(P2)} and the continuity of the map $[t \mapsto \xi(t)]$ we have that for any $t>0$ sufficiently small
    \begin{description}
      \item[bp1] if $\xi(t) < \rho_L$, then $p(\xi(t)) \equiv p(\rho_L-)$;
      \item[bp2] if $\rho_L < \xi(t)$, then $p(\xi(t)) \equiv p(\rho_L+)$.
    \end{description}
The case $\xi(t) \equiv \rho_L$ is somehow special and has to be studied separately for each specific case.

Second, when the solution is nonclassical, due to the finite speed of propagation of the waves,  the assumption~\textbf{(P2)} and properties~\textbf{bp1} and~\textbf{bp2}, we have
\begin{description}
  \item[np1] if $\mathcal{R}\left[\rho_L, \hat\rho(\bar p)\right](x) \equiv \rho_L$ for $x<0$, then $\bar p = f(\rho_L) \in \left[p(\rho_L+), p(\rho_L-)\right]$;
  \item[np2] if $\bar p \ne f(\rho_L)$ and $\rho_L < \hat\rho(\bar p)$, then $\bar p = p(\rho_L+)$;
  \item[np3] if $\bar p \ne f(\rho_L)$ and $\hat\rho(\bar p) < \rho_L$, then $\bar p = p(\rho_L-)$;
  \item[np4] if $p$ is continuous in $\rho_L$, namely $p(\rho_L-)=p(\rho_L+)$, then $\bar p = p(\rho_L)$.
\end{description}

Now we start the description of the possible cases and we proceed as follows. First, we show that for any initial datum satisfying~\textbf{(C$i$)}, $i=1,\ldots,5$, the problem actually has a unique solution and that the solution is classical. Second, we take into consideration the corresponding case~\textbf{(N$i$)}, for which we prove that the classical solution is not suitable and that there exists a unique nonclassical solution.

In general the solutions to the constrained Riemann problem~\ref{eq:constrianedRiemann} are not self--similar. All the cases listed below describe self--similar solutions because we let the solutions evolve only on a small interval of time.

    \begin{enumerate}[\textbf{(N4a)}]
      \item[\textbf{(C1)}]  In this case $\left[ (t,x) \mapsto \mathcal{R}[\rho_L, \rho_R](x/t) \right]$ performs a shock with negative speed $\sigma(\rho_L, \rho_R)$ and satisfies~\ref{eq:constrianedRiemann2} because $f(\rho_R) \le p(\rho_L+)$ and $p(\xi(t)) \equiv p(\rho_L+)$ by~\textbf{bp2}. Assume that there exists a nonclassical solution of the form~\ref{eq:nonclassicalsol}. Observe that the assumptions $\rho_L < \rho_R$ and $f(\rho_R) < f(\rho_L)$ together imply that $\bar\rho < \rho_R$.  Then $\check\rho(\bar p) \le \bar\rho < \rho_R$ and $\mathcal{R}[\check\rho(\bar p),\rho_R]$ is given by a shock with non negative speed if and only if $\bar p \le f(\rho_R)$, or equivalently, $\rho_R \le \hat\rho(\bar p)$. As a consequence, $\bar p \le f(\rho_R) < f(\rho_L)$, $\rho_L < \rho_R \le \hat\rho(\bar p)$ and by~\textbf{np2} $\bar p$ coincides with $p(\rho_L+)$. In conclusion we have $\bar p \le f(\rho_R) \le p(\rho_L+) = \bar p$, namely $f(\rho_R) =\bar p$ and the nonclassical solution coincides with the classical one.
      \item[\textbf{(N1)}] In this case $\left[ (t,x) \mapsto \mathcal{R}[\rho_L, \rho_R](x/t) \right]$ does not satisfy~\ref{eq:constrianedRiemann2} because $f(\rho_R) > p(\rho_L+)$, see case~\textbf{(C1)}. Therefore, there does not exist any classical solution and we can consider only nonclassical solutions of the form~\ref{eq:nonclassicalsol}. If $p$ is continuous in $\rho_L$, then by~\textbf{np4} we have that $\bar p = p(\rho_L)$. If $p$ experiences a jump at $\rho_L$ then, one may wonder which value in $\left[p(\rho_L+), p(\rho_L-)\right]$ has to be chosen as $\bar p$. As in the case~\textbf{(C1)}, the assumptions  imply that $\bar\rho < \rho_R$ and then that $\check\rho(\bar p) < \rho_R$ and $\bar p \le f(\rho_R)$. Then $\bar p$ is strictly smaller than $f(\rho_L)$ and $\hat\rho(\bar p) > \rho_L$. As a consequence, property~\textbf{np2} forces us to choose the unique possible value of $\bar p$, which is $p(\rho_L+)$.
      \item[\textbf{(C2)}] In this case $\left[ (t,x) \mapsto \mathcal{R}[\rho_L, \rho_R](x/t) \right]$ performs a shock with non negative speed $\sigma(\rho_L, \rho_R)$ and it satisfies~\ref{eq:constrianedRiemann2} because $f(\rho_L) \le p(\rho_L+)$. Assume that there exists a nonclassical solution of the form~\ref{eq:nonclassicalsol}. Observe that the assumptions $\rho_L < \rho_R$ and $f(\rho_R) \ge f(\rho_L)$ together imply that $\bar\rho > \rho_L$.  Then $\hat\rho(\bar p) \ge \bar\rho >\rho_L$ and $\mathcal{R}[\rho_L, \hat\rho(\bar p)]$ is given by a shock with non positive speed if and only if $\bar p \le f(\rho_L)$. Thus $\bar p \le f(\rho_L) \le p(\rho_L+)$ and this implies by~\ref{eq:nonclassicalsol3} that $\bar p = f(\rho_L) = p(\rho_L+)$ and that the nonclassical solution coincides with the classical one.
      \item[\textbf{(N2)}] In this case $\left[ (t,x) \mapsto \mathcal{R}[\rho_L, \rho_R](x/t) \right]$ does not satisfy~\ref{eq:constrianedRiemann2} because $f(\rho_L) > p(\rho_L-)$, see case~\textbf{(C2)}. Therefore, there does not exist any classical solution and we can consider only nonclassical solutions of the form~\ref{eq:nonclassicalsol}. As in the case~\textbf{(C2)},  the assumptions imply $\hat\rho(\bar p) \ge \bar\rho >\rho_L$. Furthermore, by~\ref{eq:nonclassicalsol3} we have $\bar p \le p(\rho_L-) < f(\rho_L)$, and as a consequence, property~\textbf{np2} forces us to choose $\bar p = p(\rho_L+)$.
      \item[\textbf{(C3)}] In this case $\left[ (t,x) \mapsto \mathcal{R}[\rho_L, \rho_R](x/t) \right]$ performs a possible null rarefaction on the right of the constraint and it satisfies~\ref{eq:constrianedRiemann2} because $f(\rho_L) \le p(\rho_L+)$. Assume that there exists a nonclassical solution of the form~\ref{eq:nonclassicalsol}. Since $\rho_L \le \bar\rho \le \hat\rho(\bar p)$, $\mathcal{R}[\rho_L,\hat\rho(\bar p)]$ is given by a shock that has non positive speed if and only if $\bar p \le f(\rho_L)$. Therefore $\bar p \le f(\rho_L) \le p(\rho_L+)$ and this by~\ref{eq:nonclassicalsol3} implies that $\bar p = f(\rho_L) = p(\rho_L+)$ and that the nonclassical solution coincides with the classical one.
      \item[\textbf{(N3)}] In this case $\left[ (t,x) \mapsto \mathcal{R}[\rho_L, \rho_R](x/t) \right]$ does not satisfy~\ref{eq:constrianedRiemann2} because $f(\rho_L) > p(\rho_L-)$, see case~\textbf{(C3)}. Therefore, there does not exist any classical solution and we can consider only nonclassical solutions of the form~\ref{eq:nonclassicalsol}. By hypothesis and~\ref{eq:nonclassicalsol3} we have $f(\rho_L) > p(\rho_L-) \ge \bar p$. Therefore $\rho_L \le \bar\rho < \hat\rho(\bar p)$ and by~\textbf{np2} we have $\bar p = p(\rho_L+)$.
      \item[\textbf{(C4)}] In this case $\left[ (t,x) \mapsto \mathcal{R}[\rho_L, \rho_R](x/t) \right]$ performs a rarefaction with speeds between $\lambda(\rho_L) <0$ and $\lambda(\rho_R)\ge0$ and it satisfies~\ref{eq:constrianedRiemann2} because $f(\bar\rho) = p(\rho_L+)$ implies that $p(\rho) = f(\bar\rho)$ for all $\rho \le \rho_L$. Moreover, it implies also that $p$ is continuous in $\rho_L$ and therefore, by~\textbf{np4}, any nonclassical solution of the form~\ref{eq:nonclassicalsol} must have $\bar p = p(\rho_L) = f(\bar\rho)$, but in this case the nonclassical solution coincides with the classical one.
     \item[\textbf{(N4)}] In this case $\left[ (t,x) \mapsto \mathcal{R}[\rho_L, \rho_R](x/t) \right]$ does not satisfy~\ref{eq:constrianedRiemann2} because $f(\bar\rho) > p(\rho_L-)$, see case~\textbf{(C4)}. Therefore, there does not exist any classical solution and we can consider only nonclassical solutions of the form~\ref{eq:nonclassicalsol}.
         \begin{enumerate}[\textbf{(N4a)}]
     \item[\textbf{(N4a)}] By assumption and~\ref{eq:nonclassicalsol3} $f(\rho_L) < p(\rho_L+) \le \bar p$ and therefore $\hat\rho(\bar p) < \rho_L$ and by~\textbf{np3} we have $\bar p = p(\rho_L-)$.
     \item[\textbf{(N4b)}] By assumption and~\ref{eq:nonclassicalsol3} $f(\rho_L) > p(\rho_L-) \ge \bar p$ and therefore $\hat\rho(\bar p) > \rho_L$ and by~\textbf{np2} we have $\bar p = p(\rho_L+)$.
         \end{enumerate}
     \item[\textbf{(C5)}] In this case $\left[ (t,x) \mapsto \mathcal{R}[\rho_L, \rho_R](x/t) \right]$ performs a possible null rarefaction on the left of the constraint and it satisfies~\ref{eq:constrianedRiemann2} because $f(\rho_R) \le p(\rho_L-)$ and $p(\xi(t)) \equiv p(\rho_L-)$ by~\textbf{bp1}. Assume that there exists a nonclassical solution of the form~\ref{eq:nonclassicalsol}. Since by assumption and~\ref{eq:nonclassicalsol3} $\bar p \ge p(\rho_L+) > f(\rho_L)$, we have $\hat\rho(\bar p) < \rho_L$ and by~\textbf{np3} $\bar p = p(\rho_L-)$, but in this case the nonclassical solution coincides with the classical one.
     \item[\textbf{(N5a)}] In this case $\left[ (t,x) \mapsto \mathcal{R}[\rho_L, \rho_R](x/t) \right]$ does not satisfy~\ref{eq:constrianedRiemann2} because $f(\rho_R) > p(\rho_L-)$, see case~\textbf{(C5)}. Therefore, there does not exist any classical solution and we can consider only nonclassical solutions of the form~\ref{eq:nonclassicalsol}.
    \begin{enumerate}[\textbf{(N5b)}]
    \item[\textbf{(N5b)}] By assumption and~\ref{eq:nonclassicalsol3}, $f(\rho_L) < p(\rho_L+)\le \bar p$ and therefore $\hat\rho(\bar p) < \rho_L$ and by~\textbf{np3} we have $\bar p = p(\rho_L-)$.
     \item[\textbf{(N5b)}] By assumption and~\ref{eq:nonclassicalsol3}, $f(\rho_L) > p(\rho_L-) \ge \bar p$ and therefore $\hat\rho(\bar p) < \rho_L$ and by~\textbf{np2} we have $\bar p = p(\rho_L+)$.
    \end{enumerate}
    \end{enumerate}

\subsubsection{Cases in which uniqueness is violated}

Now we list the ``pathological'' cases, where we have more than one admissible solution. We stress once again that a necessary condition for non--uniqueness is  $p(\rho_L-) \ne p(\rho_L+)$ and $p(\rho_L-) \ge f(\rho_L) \ge p(\rho_L+)$.
\begin{description}
\item[(CN2)] If $\rho_L < \rho_R$, $f(\rho_L) \le f(\rho_R)$ and $p(\rho_L+)< f(\rho_L) \le p(\rho_L-)$, then the classical solution $\left[ (t,x) \mapsto \mathcal{R}[\rho_L, \rho_R](x/t) \right]$, which consists of a shock with non negative speed, as well as the nonclassical solution~\ref{eq:nonclassicalsol}, with $\bar p = p(\rho_L+)$, are distinct solutions of~\ref{eq:constrianedRiemann}.
\item[(CN3)] If $\rho_R \le \rho_L\le \bar\rho$ and  $p(\rho_L+)< f(\rho_L) \le p(\rho_L-)$, then the classical solution $\left[ (t,x) \mapsto \mathcal{R}[\rho_L, \rho_R](x/t) \right]$, which consists of a possible null rarefaction on the right of the constraint, as well as the nonclassical solution~\ref{eq:nonclassicalsol}, with $\bar p = p(\rho_L+)$, are distinct solutions of~\ref{eq:constrianedRiemann}.
\item[(NNN4)] If $\rho_R \le \bar\rho < \rho_L$, $p(\rho_L-) \ne p(\rho_L+)$ and $p(\rho_L+) \le f(\rho_L) \le p(\rho_L-)$, then the nonclassical solutions of the form~\ref{eq:nonclassicalsol} which corresponds to $\bar p$ in the set $ \{p(\rho_L+), f(\rho_L),  p(\rho_L-)\}$ satisfy~\ref{eq:constrianedRiemann}. This is the situation considered in the Example~2 in \cite{AndreianovDonadelloRosini}. Observe that such solutions are distinct as far as they correspond to distinct constraint levels $\bar p$, and that in any case there exist at least two distinct nonclassical solutions.
\item[(CNN5)] If $\bar\rho < \rho_R \le \rho_L$, $f(\rho_R) \le p(\rho_L-)$, $p(\rho_L-) \ne p(\rho_L+)$ and $p(\rho_L+) \le f(\rho_L)$, then the classical solution $\left[ (t,x) \mapsto \mathcal{R}[\rho_L, \rho_R](x/t) \right]$, which consists of a possible null rarefaction on the left of the constraint, as well as the nonclassical solutions of the form~\ref{eq:nonclassicalsol} corresponding to $\bar p \in \{p(\rho_L+), f(\rho_L)\}$ satisfy~\ref{eq:constrianedRiemann}. Observe that the two nonclassical solutions are distinct as far as they correspond to distinct constraint levels $\bar p$, and that in any case there exist at least two distinct solutions, one classical and one nonclassical.
\item[(NNN5)] If $\bar\rho < \rho_R < \rho_L$, $f(\rho_R) > p(\rho_L-) \ge f(\rho_L) \ge p(\rho_L+)$ and $p(\rho_L-) \ne p(\rho_L+)$, then the nonclassical solutions of the form~\ref{eq:nonclassicalsol} corresponding to $\bar p \in \{p(\rho_L+), f(\rho_L),  p(\rho_L-)\}$ satisfy~\ref{eq:constrianedRiemann}. Observe that such solutions are distinct as far as they correspond to distinct constraint levels $\bar p$, and that in any case there exist at least two distinct nonclassical solutions.
\end{description}

\subsection{Proof of Proposition~\ref{prop:Riemann}}\label{sec:technicalRiemann}

\begin{enumerate}[\textbf{(R1)}]
\item[\textbf{(R1)}] Any solution given by $\mathcal{R}^\star$ coincides on each side of the constraint with a solution given by the classical Riemann solver $\mathcal{R}$. Therefore it satisfies the Rankine--Hugoniot jump condition along any of its discontinuities away from the constraint. Finally, by definition of $\hat\rho$ and $\check\rho$, it satisfies the Rankine--Hugoniot jump condition also along the constraint.
\item[\textbf{(R2)}] It is clear by the proof of Proposition~\ref{prop:riemann}.
\item[\textbf{(R3)}] It is proved as in~\textbf{(R1)} since any classical solution is in $\BV$.
\item[\textbf{(R4)}] As was proved in~\cite{ColomboGoatinConstraint}, $\mathcal{R}^\star $ is continuous on $\mathcal{C}\cup\mathcal{N}$. If $(\rho_L,\rho_R)$ is not in $\mathcal{C}\cup\mathcal{N}$ then $ p$ experiences a jump at $\xi=\rho_L$. Therefore, the local in time solutions of the Riemann problem for the initial conditions $(\rho_L+\varepsilon,\rho_R)$ and $(\rho_L-\varepsilon,\rho_R)$ are different and only one of the two converges to  $\mathcal{R}^\star[\rho_L, \rho_R]$ as $\varepsilon>0$ goes to zero.
\item[\textbf{(R5)}] We first stress once again that we can discuss the consistency property of our Riemann solvers only locally in time because, in general, the solutions may be not even self--similar globally in time. However, locally in time, the efficiency of the exit can be assumed to be constant and it is thus sufficient to proceed as in~\cite{ColomboGoatinConstraint}.
\item[\textbf{(R6)}] It is clear by the proof of Proposition~\ref{prop:riemann}.
\end{enumerate}

\section{Aknowledgements}
All the authors are supported by French ANR JCJC grant CoToCoLa and Polonium 2014 (French-Polish cooperation program) No.331460NC. The second author is also supported by the Universit\'e de Frache-Comt\'e, soutien aux EC 2014. The last author is also supported by ICM, University of Warsaw, and by Narodowe Centrum Nauki, grant 4140.

\bibliographystyle{plain}
    \bibliography{bibliography}

\end{document}